\magnification=\magstep1
\dimen100=\hsize
\parskip= 6pt

\font\ninerm=cmr10 at 10truept
\font\eightrm=cmr8
\font\sc=cmcsc10

\font\tenmsy=msbm10
\font\sevenmsy=msbm7
\font\fivemsy=msbm5
\newfam\msyfam
\def\msy{\fam\msyfam\tenmsy}
\textfont\msyfam=\tenmsy
\scriptfont\msyfam=\sevenmsy
\scriptscriptfont\msyfam=\fivemsy

\def\bbc{{\msy C}}

\def\bbi{{\msy I}}

\def\bbp{{\msy P}}

\def\bbz{{\msy Z}}
\def\za{\vrule height6pt width4pt depth1pt}

\font\aa=eufm10

\def\Got#1{\hbox{\aa#1}}

\def\Ra*{(R_{a})_{*}}
\def\ada-1{ad_{a^{-1}}~}
\def\p-1u{\pi^{-1}(U)}
\def\gg{{\Got g}}

\def\gs{{\Got s}}

\def\gl{{\Got l}}

\def\gsl{{\Got sl}}

\font\svtnrm=cmr17

\centerline{\svtnrm Surfaces and the Sklyanin bracket} 
\bigskip
\centerline{\sc  J.C. Hurtubise and E. Markman}
\footnote{}{\ninerm The  first  author of this article would like to thank 
NSERC and FCAR for their support}
\footnote{}{\ninerm The second author was partially supported by NSF grant 
number DMS-9802532}
\bigskip
\centerline{\vbox{\hsize = 5.85truein
\baselineskip = 12.5truept
\eightrm
{\sc Abstract}: We discuss the Lie Poisson groups structures associated to
splittings of the loop group $LGL(N,\bbc)$, due to  Sklyanin. 
Concentrating on the finite dimensional leaves of the associated Poisson 
structure, we show that the geometry of the leaves is
intimately related to a complex algebraic ruled surface with
a $\bbc^\times$-invariant Poisson structure.  In particular, Sklyanin's
Lie Poisson structure admits a suitable abelianisation, once one passes to
an appropriate spectral curve. The Sklyanin structure is then equivalent 
to one considered by Mukai, Tyurin and Bottacin on a moduli space 
of sheaves on the Poisson surface. 
The abelianization procedure gives rise to
natural Darboux coordinates for these leaves, as well as separation 
of variables for  the  integrable Hamiltonian systems 
 associated to invariant functions on the group. }}

{\bf 1. Introduction}


The aim of this note is to close up a circle of ideas linking
algebraically integrable systems associated to 
loop groups and loop algebras on $Gl(N, \bbc)$ or $Sl(N,\bbc)$,
to symmetric products of 
certain symplectic surfaces (more properly, Hilbert schemes of zero
dimensional subschemes  on
the surfaces). The general idea is to show that the phase spaces
of  these systems 
are birationnally symplectomorphic  to the Hilbert schemes, in such a way
that the
leaves of the Lagrangian foliation are given by the  space of 
divisors on {\it spectral curves}; these spectral curves  lie in the
surface, and the
inclusion of Lagrangian leaves (symmetric products of the curves)
into the symplectic leaves (symmetric product of the surface) is induced
by the 
inclusion of the curves into the surface. This isomorphism, when made
explicit, 
gives   simple separations of variables for the systems.

Integrable systems  on loop algebras and loop groups include as special
cases  
most of the frequently studied integrable systems. Indeed, appropriate
choices of rank,
location of poles, residues at the poles, etc, give us most 
of the classical systems: the Neumann oscillator, the various tops, the
finite gap solutions to the KdV, the NLS, the CNLS and the Boussinesq
 equations, the various Gaudin models, the Landau-Lifshitz equation are 
 just some of the examples.
References include  the book [FT], the survey [RS2], and the references
therein,  or the articles [Mo, AvM, RS1, AHP, 
HaH].
All of these systems are associated to splittings of the loop algebra
$L\gg= L\gg_+ \oplus
L\gg_-$ or the corresponding local decompositions of the loop groups 
$LG=LG_+LG_-$; this splitting gets encoded in terms of  $r$-matrices. One
has three main 
types of splitting, given by the rational, trigonometric and elliptic 
$r$-matrices. The splittings allow us to define a bracket on the loop
algebra
(linear, or Lie Poisson bracket), and on the group (quadratic, or Sklyanin 
bracket [Sk1, S]). For all of these brackets, there are integrable systems
whose 
Hamiltonians are the coefficients of the equations defining the spectral
curve of the loop.

The linear brackets admit an important generalisation: the generalised
Hitchin systems  [Hi1, Hi2, Bo1, M]. These systems are defined
on moduli spaces of pairs
$(E,\phi)$, where $E$ is a principal $G$ bundle over a Riemann surface
$\Sigma$,
and $\phi$ is a meromorphic 1-form valued section of the associated
adjoint bundle. The 
$r$-matrix systems correspond to the cases when the bundles are rigid: 
choosing trivialisations,  the
sections $\phi$ take their values in suitable subspaces of the loop
algebra.
For $Sl(N,\bbc)$, these cases occur when the curve $\Sigma$ is rational
(rational $r$-matrix),
elliptic (elliptic $r$-matrix) or rational nodal (trigonometric
$r$-matrix).
In short, there is a table: (L.P.: Lie Poisson, Sk.: Sklyanin)

\bigskip
\vbox {\offinterlineskip
\hrule
\halign{&\vrule#&\strut\ #\hfil\ \cr
height2pt&\omit&&\omit&&\omit&&\omit&\cr
&&&Rational curve& & Rat'l nodal curve & &Elliptic curve & & General
curve& \cr
height2pt&\omit&&\omit&&\omit&&\omit&\cr
\noalign{\hrule}
height2pt&\omit&&\omit&&\omit&&\omit&\cr
& Algebra & &Rational L.P. & &Trigonometric L.P.& & Elliptic L.P.& &
Generalised Hitchin& \cr
& Group& &Rational Sk.& &Trigonometric Sk.&  & Elliptic Sk.& & - & \cr
height2pt&\omit&\cr}
\hrule }

The first (general) case of the isomorphism of these systems with Hilbert
schemes  was given in [AHH] (special cases were considered in  [NV]).  
The case considered in [AHH] 
is the linear Lie-Poisson (rational r-matrix) bracket on the dual
$L\gg_+^*$ of the loop algebra
$L\gg_+$ of polynomial loops in the Lie algebra $gl(N,\bbc)$. 
The (reduced) symplectic leaves are 
reductions of coadjoint orbits in $L\gg_+^*$ by the 
action of $gl(N,\bbc)$.
It was shown that one has natural Darboux
coordinates on the symplectic leaves, which establish a
birational symplectomorphism between the 
symplectic leaves  and the Hilbert scheme of a rational symplectic 
surface.  The surface is a blow-up of the total space of the line bundle
$K_{\bbp_1}(D)$ over $\bbp_1(\bbc)$, where
$D$ is the divisor of poles of the meromorphic section. The blow up is
taken at the  intersection of the spectral curve with the inverse image of 
the divisor $D$. 
This picture was generalised in [HK], following [H], to cover the case of
the generalised 
Hitchin systems, which then specialises to the case of linear elliptic and 
linear trigonometric $r$-matrices. The surfaces one obtains are blow-ups
of line bundles over the base curve $\Sigma$.

For the quadratic brackets, the rational case was treated by Scott in
[Sc]; one  again obtains a blow up of a line bundle over $\bbp_1(\bbc)$, 
but now blowing up at the intersection of the spectral curve with 
the zero section. 
The main purpose of this note is to prove the corresponding
result in the  elliptic and trigonometric cases, which completes the list
of six cases enumerated above; at the same time, we give a general,
simplified exposition applicable to all of the cases, as well as a 
geometric 
indication of why the quadratic brackets do not generalise to higher 
genus curves. 
The finite dimensional leaves of the quadratic systems are shown to
be symplectic leaves of moduli spaces of Higgs bundles with a rigid
$Sl(N,\bbc)$-bundle and a generically invertible Higgs field. 
The abelianization of the moduli space of Higgs bundles
amounts to its realization as a moduli space of sheaves on a ruled
surface. The Poisson structure is then equivalent to one considered by 
Mukai, Tyurin and Bottacin. In turn, this structure is induced by a 
structure on the ruled surface. A similar situation holds for the 
 generalized Hitchin system. The structure there is invariant under 
 translation on the fibers of the ruling; this translates into  a linear 
 Poisson structure on
the space of matrices. In contrast, the Sklyanin brackets, which  
are  quadratic  on the space of matrices, correspond 
to a $\bbc^\times$-invariant Poisson structure on the surface 
(see Theorem 3.24 and the following remark). 
The ruled  surface admits such a $\bbc^\times$-invariant  Poisson 
structure when the base curve is rational or elliptic.
Indeed, one of the virtues of this approach is that it classifies to a
certain degree the possible Poisson brackets on the loop algebras or 
loop groups: the linear and quadratic brackets, between them, 
exhaust the possible Poisson structures on the ruled surface, 
when the base curve is elliptic or rational. 
When the genus of the base curve is higher, one only has Poisson 
structures on the surface which correspond to  the
linear brackets, so that the table above arises in some sense as a
consequence of the classification of Poisson structures on a ruled
surface.
Indeed, we will see that, in trying to generalise the Sklyanin
structure to the higher genus base curve case, 
we do not obtain integrable systems on symplectic varieties, 
but rather  on varieties equipped with a degenerate closed form. 

We will generalize the construction to arbitrary reductive groups
in a separate paper [HM]. 
Rigid bundles on an elliptic curve do not exist 
for simple groups of type other than $A_n$ ([FM] Theorem 5.13). 
This translates  into the fact that integrable systems with a quadratic 
Poisson structure can not be constructed on loop groups of type
other than $A_n$ [BD]. 
Nevertheless, allowing the principal bundle to deform, we do
get integrable systems with a quadratic bracket on moduli spaces of
pairs $(P,g)$ consisting of a principal $G$-bundle $P$ and 
a meromorphic section of its adjoint group bundle. These are equivalent 
to those obtained from the dynamic r-matrix formalism [F, EV].

 The paper is organized as follows: 
Section 2 recalls some basic facts about the Sklyanin
bracket, specialised to the context of loop algebras. 
In section 3 we exhibit the birational isomorphism between 
the Sklyanin systems and Hilbert schemes of symplectic surfaces, in the 
elliptic case. We show that the Sklyanin systems are supported
on the moduli space of Higgs pairs and the Poisson structure 
is a natural translation of the Mukai-Tyurin 
Poisson structure, via the one-to-one correspondence between
a Higgs pair and its spectral data. 
Sections 4 and 5 discuss the trigonometric and rational cases, 
respectively.
Section 6 discusses the case of general base curves.

{\it Acknowledgments:} 
We  would like to thank Pavel Etingof, Robert Friedman, John Harnad and
Alexander Polishchuk for helpful communications.

{\bf 2. Loop groups and the Sklyanin bracket}

In this section we recall certain basic facts about the Sklyanin bracket
in the 
case which concerns us here, that of a loop group with a (local)
decomposition into 
a product of two subgroups. Our main aim is to identify the finite 
dimensional symplectic leaves. After normalisatrion, these correspond to 
meromorphic maps on the curve.  A good general reference is [RS2], section
II.12

Let $G$ be a reductive complex group, thought of as a subgroup of some
$Gl(N, \bbc)$,
and consider the loop group $LG$ of analytic mappings of the circle into
$G$. 
Let us now assume that this circle is embedded into a Riemann surface
$\Sigma$, and bounds a 
disk $U_-$ whose centre $p$ will correspond to $z= \infty$; set $U_+=
\Sigma\backslash p$. We furthermore choose a transition
matrix $T(z)$ defined over the punctured disk, possibly taking values  
 in a larger group $\hat G$, but such that  
$G$ is a normal subgroup of $\hat G$.
We assume that the Lie algebra 
$L\gg$ of $LG$ splits into a sum 
$L\gg = L\gg_- \oplus L\gg_+$, where 
\item{-} $L\gg_-$ consists of holomorphic functions defined over $U_-$
into $\gg$,
\item{-} $L\gg_+$ consists of elements of the form
 $Ad_{T (z)} a(z)$
 where $a(z)$  is a holomorphic function from $U_+$  into $\gg$. 
 
If $P$ is the principal $\hat G$-bundle defined by $T$, 
the existence of a splitting is equivalent to the vanishing of
the first cohomology of the adjoint bundle 
$H^1(\Sigma, P(\gg)) = 0$. When $\hat G = G$,  
standard deformation theory shows that the vanishing is tantamount to $P$,
defined by the transition matrix $T$, being rigid. Rigidity
of $P$ means that any small deformation of $P$ is isomorphic to $P$. 
If $\hat \gg$ is reductive with  semi-simple part $\gg$, then 
the vanishing of $H^1(\Sigma, P(\gg))$ means that the associated
principal $G_{\rm adj}$-bundle is rigid, where $G_{\rm adj}$
is the adjoint group. In other words, the 
deformations of $P$ all arise from the center of $\hat G$.
Such rigidity  happens only for curves of genus zero or one: 
if the curve is rational, $G$ can be any reductive group; 
if the curve is elliptic, $G$ must be of type $A_n$ ([FM] Theorem 5.13).
In the elliptic case we choose $G=Sl(n,\bbc)$ and $\hat G = Gl(n,\bbc)$. 
Uniqueness in the decomposition is linked 
to the absence of sections of the adjoint bundle: 
$H^0(\Sigma,P(\gg)) =0$.

We write the corresponding (local) decomposition of groups as 
$$LG = LG_+\cdot LG_-.\eqno (2.1)$$
and decompose elements of $LG$ as follows:
$$ g= g_+\cdot g_-^{-1}.\eqno (2.2)$$
 Let $\omega$ be a non-vanishing holomorphic one-form on $U_-$. If
$(a,b) \mapsto
tr(ab)$ 
denotes the Killing form on $\gg$, we can define a pairing on $L\gg$ by 
$$<f,g> = \oint tr(fg)\omega.\eqno(2.3)$$
With respect to this pairing, both $L\gg_+$ and $L\gg_-$ 
are isotropic.

Now let $P_+,P_-$ denote the projections of $L\gg$ onto $L\gg_+$ and
$L\gg_-$
respectively. We set
$$R= P_+-P_-.\eqno (2.4)$$
If $\psi$ is a function on $LG$,   the left derivative $D\psi$ and the
right derivative $D'\psi$ in $L\gg$ are defined 
at $g$ by 
$$\eqalign{<D\psi (g) , h> &= {d\over dt} \psi( exp(th)\cdot g)|_{t=0},
\ \ \ \forall h \in L\gg,\cr
<D'\psi (g) , h> &= {d\over dt} \psi(g\cdot exp(th))|_{t=0},
\ \ \ \forall h \in L\gg.}\eqno (2.5)$$
In terms of the Maurer-Cartan forms $\theta = dg\cdot g^{-1},
\theta' = g^{-1}dg$, thought of as maps $TLG\rightarrow L\gg$, we
have $$\eqalign {D\psi &= (\theta^{-1})^* (d\psi),\cr
D'\psi &= (\theta'^{-1})^* (d\psi).}\eqno (2.6)$$
At an element $g$ of $LG$, $$D\psi= Ad_g(D'\psi).\eqno (2.7)$$
The Sklyanin bracket of $\psi$ and $\phi$ is then defined by
$$\{\psi, \phi\} = {1\over 2}<R(D\psi), D\phi> - {1\over
2}<R(D'\psi),D'\phi>. \eqno (2.8)$$
Alternately, one can write
$$\eqalign{ \{\psi,\phi\} =&  \ {1\over 2}< D\psi_+- D\psi_--
 Ad_g(D'\psi_+- D'\psi_-), D\phi>\cr
=& \ < D\psi_+ -  Ad_g(D'\psi_+ ), D\phi>\cr
=& \ < -  D\psi_-+ Ad_g(  D'\psi_-), D\phi>,\cr
}\eqno (2.9)$$

\noindent{\it  The dressing action.}
One has a (right) action of $LG_-\times LG_+$  on $LG$, defined by
$$ g(h_+, h_-)  = ((gh_+h_-^{-1}g^{-1})_+)^{-1}gh_+=
((gh_+h_-^{-1}g^{-1})_-)^{-1}gh_- \eqno(2.10) $$
For $g$ lying in $LG_-$, for example, the action of $h_-$ is trivial, and
the action of
 $h_+$ is  
$$g\mapsto ((gh_+)_-)^{-1} = ((gh_+)_+)^{-1}gh_+.$$
Referring to (2.2), this is in essence the right action of $h_+$
on $g$, followed by projection to $G_-$. A more conceptual definition
can be given of 
the dressing action, 
involving projection of simple flows on a larger space; 
see, e.g., proposition 12.20 of [RS2].
Set, for an element $\zeta$ of $L\gg$:
$$\zeta_+ = P_+(\zeta),\quad \zeta_- = P_-(\zeta).$$
We note that this differs from the infinitesimal version of (2.2) by  a 
sign.
If $\xi_+\in L\gg_+,
\xi_-\in L\gg_-$, the actions of the one parameter subgroups
$exp (\epsilon\xi_\pm)$ are given by 
$$g\mapsto (1+\epsilon (Ad_g\xi_+)_- + O(\epsilon^2))g,\quad
g\mapsto (1+\epsilon (Ad_g\xi_-)_+ + O(\epsilon^2))g.
\eqno(2.11)$$
In other words, for a function $f$, if $v_{\xi_\pm}$ denotes the  
vector field corresponding to $\xi_\pm$,
$$v_{\xi_+}(f) = <  (Ad_g(\xi_+))_-, Df>, \quad
v_{\xi_-}(f) = <  (Ad_g(\xi_-))_+, Df>.\eqno(2.12)$$
One has (see [RS2], [S]):

{\sc Theorem} (2.13) \tensl The symplectic leaves of the Sklyanin bracket 
are given by the orbits of the dressing action. \tenrm


 From now on, we take  $G= Gl(n,\bbc)$ or $Sl(n,\bbc)$. We can analyse the
finite dimensional 
symplectic leaves as follows:

{\sc Theorem (2.14) } \tensl
\item {a)} The finite dimensional leaves in $LGL(n,\bbc)$
are orbits of elements of the form $f(z) g(z) $
 where $f(z)$
is a scalar function, $g(z)$ has a pole of finite order at p 
and  $T(z)g(z)T(z)^{-1}$ is a meromorphic
 matrix-valued function on $\Sigma-\{p\}$, with a finite number of poles.

\item {b)} The location and the order of the poles is constant along the
orbit.

\item {c)} The points over which $det (g) $ vanishes are constant along
the 
orbit.
\tenrm

{\sc Proof:} Normalise one of the matrix coeficients of $g^{-1}$, say
$g^{-1}_{11}$, to 1; this 
accounts for the function $f$. 
Now look at the action of $L\gg_-$, to analyse the pole at $\infty$. 
Filter $L\gg_-$ as $...(L\gg_-)_{-n-1}\subset (L\gg_-)_{-n}\subset ....$
by order of vanishing
at $\infty$.
 Finite codimension of the stabiliser $W_-$ implies that
$(L\gg_-)_{-n}\subset W_- $ for some $n$:
there can only be a finite set of $n$'s such that the map
$(L\gg_-)_{-n}\cap W_- \rightarrow \gg=
(L\gg_-)_{-n}/(L\gg_-)_{-n-1}$ is not surjective. Next consider elements
of $(L\gg_-)_n$ of the form $z^n \epsilon_{j1}$, 
where $\epsilon_{j1}$ is the $(j,1)$-th  elementary matrix
and using the fact that $g_{ij} z^n = g_{ij} z^n \epsilon_{j1}g^{-1}_{11}$
is holomorphic 
at $\infty$ tells us that $g$ has a pole of order at most $n$.

Next we consider the action of $L\gg_+$. The stabiliser $W_+$ is an ${\cal
O}(\Sigma-\infty)$-module
and so we can localise over the points $q$ of $\Sigma-\{\infty\}$, and
consider the quotients
$(L\gg_+)_q/(W_+)_q$. Again, the finiteness of the codimension of $W_+$
tells us that only a 
finite number of these quotients are non-zero, and that at each point the
quotient module is supported on 
some finite formal neighbourhood of the point. This tells us that there is
a function $f$ on
$\Sigma-\{\infty\}$ with finite order poles at a finite number of points 
such that $f\cdot L\gg_+$ maps to zero in
$(L\gg_+)_q/(W_+)_q$ 
at each point $q$ of 
$\Sigma-\{\infty\}$, and so $f\cdot L\gg_+\subset W_+$. Again this tells
us that 
the components $f g_{ij}$ lie in ${\cal O}(\Sigma-\{\infty\})$. 

For b), note that we can write the dressing action of $\xi_+$, using
(2.11), as:
$$\dot g = (g\xi_+g^{-1})_-g =  g\xi_+ - (g\xi_+g^{-1})_+g,$$
and similarily, for $\xi_-$
$$\dot g = (g\xi_-g^{-1})_+g =  g\xi_- - (g\xi_-g^{-1})_-g,$$
from which it follows that the poles of $\dot g$ are included amongst
those 
of $g$: the dressing action preserves 
the singularities of $g$.  Part c) follows by considering the explicit 
form (2.10)
of the dressing action.
 \hfill\za

{\sc Remark} In keeping with our interpretation of the decomposition
(2.1) in terms of
 a holomorphic bundle, we note that there 
are two operations we can perform on the brackets:

\noindent(2.15) {\it Changing trivialisations.} As we are thinking of our 
splitting 
$L\gg= L\gg_+\oplus L\gg_-$ in terms of sections of a rigid bundle $E$, we
should take advantage of this 
and allow ourselves to change trivialisations. Let us then set $\hat T =
T_-T_+^{-1}$, where 
$T_\pm:U_\pm\rightarrow G$ are holomorphic maps. Define $\widehat
{L\gg}_+$ to be 
$T_-(L\gg_+)T_-^{-1}$,
and define the map 
$$\eqalign{ \rho: L\gg_+\times L\gg_- &\rightarrow \widehat{L\gg}_+\times
L\gg_-\cr
(a,b)&\mapsto (T_-aT_-^{-1}, T_-bT_-^{-1}).}$$
If we define modified projections by
$$\hat P_\pm(a) = T_-P_\pm(T_-^{-1}aT_-) T_-^{-1},$$
and set $\hat R = \hat P_+- \hat P_-$, we have that the map $\rho$
intertwines the two Poisson brackets 
defined by $R, \hat R$. In other words, we may work with the
trivialisation we wish.
 
\noindent(2.16) {\it Adding points.} A given bundle of course admits not 
only several
trivialisations with respect to 
a fixed covering by open sets, but also trivialisations with respect to
different coverings.
In particular, let us suppose that we have not only a covering by 
$U_+= \Sigma-\{p\}, U_-=$ disk containing $p$,
but also a cover by $\tilde U_+= \Sigma-\{p, q\}, \tilde U_{-,p}=U_-,
\tilde U_{-,q}=$ disk containing $q$, where the two disks 
$\tilde U_{-,p},
\tilde U_{-,q}$ do not overlap. For this
 second cover, the functions on the overlaps $\tilde U_+\cap (\tilde
U_{-,p}\cup\tilde U_{-,q})$
correspond to  a sum of two copies of the loop algebra $L\gg \oplus L\gg$.
We keep the transition function $T$ on 
$\tilde U_+\cap  \tilde U_{-,p}$, and take the identity as transition
function on $\tilde U_+\cap  \tilde U_{-,q}$.
We can decompose $L\gg \oplus L\gg$ into a sum $\widetilde {L\gg}_+ \oplus
\widetilde {L\gg}_-$, where 
$\widetilde {L\gg}_+$ corresponds to sections of $ad(E)$ over $\tilde
U_+$, and $\widetilde {L\gg}_- = 
{L\gg}_{-,p}\oplus {L\gg}_{-,q}$ consists of sections over the two open
disks. There are corresponding projections 
$\tilde P_+, \tilde P_-$, and a corresponding Sklyanin bracket.
One can show that  the projection $\pi: L\gg \oplus L\gg \rightarrow L\gg$ 
onto the first factor is a Poisson map. More generally, we can add and
subtract points, which shows that
 the intrinsic object we are considering is really the space of sections
of a rigid bundle.
\bigskip

{\bf 3. Spectral curves and Abelianisation: the elliptic case}

As we have seen, we are in essence considering sections over a punctured 
disk of the automorphisms of  rigid bundle 
over a Riemann surface; the finite dimensional symplectic leaves are 
those of meromorphic sections of the automorphisms over the whole surface.
There are three cases that one can consider, those of an elliptic curve, 
a rational nodal curve, and a rational curve. This section is devoted to 
the elliptic case.

{\it 3.a  A rigid bundle on an elliptic curve.}
Let
$\Sigma $ be an elliptic
curve and $D$   a positive divisor on $\Sigma$. We will 
take  as vector bundle $E$  a stable vector bundle of rank $N$,
degree $1$, and we fix the top exterior power of $E$. This makes 
 the bundle rigid, and in fact determines the bundle. The bundle $E$ can
be defined as follows.
Let $q=  {\rm exp} (2\pi i/N)$, and set
$$I_1 ={\rm diag} (1,q, q^2,...,q^{N-1}),\quad\quad
I_2= \pmatrix {0&1&0&\dots&0\cr 0&0&1&\dots&0\cr.&.&.& &.\cr.&.&.& &.\cr
0&0&0&\dots&1\cr1&0&0&\dots&0}.\eqno(3.1)$$
Now let us represent the elliptic curve $\Sigma$ as
$\bbc/\bbz\omega_1+\bbz\omega_2$.
Puncture the curve at a point $p$. One can lift  $E$ to $\bbc$;
sections of $E$ over $\Sigma$ will be given by $N-$tuples $F$ of
$N$-valued functions
defined over the inverse image in $\bbc$ of $\Sigma-p$, satisfying:
\item {-} $F(z+\omega_i) = I_i F(z)$,
\item {-} $F$ is of the form $z^{-1/N}\cdot$holomorphic, near the inverse
images in 
$\bbc$ of the puncture $p$, where $z=0$ corresponds to $p$.

In a similar vein, sections of $End(E)$ are given by holomorphic
matrix-valued
functions $M$ (this time single-valued) on $\bbc$, satisfying
$M(z+\omega_i) =
I_iM(z)I_i^{-1}$. 
 We will consider  the  subspace 
$H^0(\Sigma, End(E)(D))$ of
 endomorphisms of $E$, with  the order of the poles bounded by the divisor 
$D$. 

We note that, because the degree is $1$, 
we have a bundle with structure group 
$Gl(N, \bbc)$. The group that we consider for our splitting, however, 
is $Sl(N,\bbc)$, with Lie algebra
$\gg = \gs\gl(n,\bbc)$, or, alternately, the group $PGL(N,\bbc)$.
This will correspond to the traceless
endomorphisms
 $End^{~0}(E)$.  We  have
$$H^0(\Sigma, End^{~0}(E))= H^1(\Sigma, End^{~0}(E))= 0,$$
so that we are indeed in the case of a unique decomposition, as in (2.1).
The decompositions we consider will thus be of the group of sections of 
$Sl(E)$ over the punctured disk, or more generally, sections of $Gl(E)$ 
with fixed determinant.
While  the bundle has been defined using automorphy factors, 
rather than transition matrices, a change of trivialisations, as 
in remark (2.15), allows one to go from one formalism to the other.

{\it 3.b The Mukai structure.} We have already seen in the previous 
sections, that the space of sections of $H^0(\Sigma, End(E)(D))$ with a 
fixed determinant is a union of symplectic leaves for the Sklyanin 
structure. We will now construct on  $H^0(\Sigma, End(E)(D))$ another 
 Poisson structure, whose 
symplectic leaves will again  be subvarieties of $H^0(\Sigma, End(E)(D))$ 
with fixed determinant. 
We will proceed by   reduction by a $\bbc^*$-action of 
a larger space ${\cal M}$ of pairs $(E',g)$ where $E' = E \otimes L$ for 
a line bundle $L$,
and $g\in  H^0(\Sigma, End(E)(D))$ is generically invertible. 
Symplectic leaves of ${\cal M}$ are determined by the zero divisor of 
the determinant of $g$, so that ${\rm det}(g)$ is fixed only up to
a scalar factor. The $\bbc^*$-action is defined by 
$$c(E,g) = (E, c\cdot g),   \eqno (3.2)$$ 
and so, up to a finite cover corresponding to action by roots of unity, 
taking the quotient by $\bbc^*$ corresponds to fixing the determinant.

We have not mentioned yet the Hamiltonians that will define our integrable
 sustems on the finite dimensional symplectic leaves; this system is 
 closely tied to the Mukai structure. The  Hamiltonians are given by
the coefficients of the defining 
equation $F$ of the spectral curve $S$ of $g\in H^0(\Sigma, End(E)(D))$: 
$$ F(z,\lambda) = det(g(z)-\lambda\bbi) = 0.\eqno (3.3)$$
In short, the Lagrangian leaves of the integrable system are given by
fixing the spectral curve.
If $D$ is the divisor of poles of $g$, the equation (3.3) defines a
compact curve $S$ embedded in the 
total space ${\cal T}$ of the line bundle ${\cal O}(D)$ over $\Sigma$;
there is an $n$-sheeted projection
$\pi:S\rightarrow\Sigma$.
One can also define a sheaf $L$ supported over the spectral curve as a
cokernel of $g-\lambda\bbi$; generically it is a line bundle over $S$,
$$0\longrightarrow \pi^*E\otimes {\cal
O}(-D){\buildrel{g-\lambda\bbi}\over
{\longrightarrow}}\pi^* E\rightarrow L
\rightarrow 0 .\eqno(3.4)$$
We have:
     
{\sc Proposition} (3.5) [H]\tensl
\item {a)} The push-down $\pi_*L$ is isomorphic to $E$. 
\item {b)} The map $g$, up to  conjugation by the global automorphisms of 
$E$, is the
push-down of 
the action on $L$ given by multiplication by the 
fiber coordinate $\lambda$.\tenrm

The automorphisms of our family of $E$'s are multiples of the identity,
so that one recovers the pair $(E,g)$ from $(S,L)$.

{\sc Proposition} (3.6) \tensl
Let ${\cal S}$ be the family of smooth curves $S'$ in the linear system of
$S$ on the surface ${\cal T}$.
Then ${\cal M}$ contains the Jacobian fibration of ${\cal S}$
(of degree $g$ line-bundles) as a Zariski open subset.  \tenrm


Let us consider the deformation theory of $L$, first as a line bundle
supported over a smooth
curve. The normal bundle of the spectral curve is given by the twist
$K_S(D)$ of 
the canonical bundle of $S$
by $\pi^*({\cal O}(D))$ and so the    space of infinitesimal deformations
of the curve is then  
$ H^0 (S, K_S(D))$. If one constrains the sections of the normal bundle to 
vanish
on the zero-section in ${\cal T}$, one gets a space of sections isomorphic
to $H^0(S, K_S)$. 
 Deformations of the line bundle, 
fixing the curve, are given by $ H^1(S,{\cal O})$. On the other hand (and 
more generally), one
can think of $L$ as a sheaf on ${\cal T}$:
deformations of $L$ as a sheaf on  ${\cal T}$ include both deformations of
its support,  and deformations of the line bundle. These
deformations are classified by the extension group $Ext^1_{{\cal
T}}(L,L)$.  
On ${\cal M}$, we have an exact sequence for the tangent bundle, 
linked to the fact that it is  the  Jacobian fibration:
$$0\rightarrow H^1(S,{\cal O})\rightarrow T{\cal M}= Ext^1(L,L)
\rightarrow H^0 (S, K_S(D))\rightarrow 0\eqno (3.7)$$

We will show that the Sklyanin structure is equivalent to one defined by 
Tyurin and Bottacin [Bo2,T] for  sheaves on a Poisson surface 
(generalizing the work of Mukai [Mu]). 
The surface that we are considering is ${\cal T}$; 
the top exterior power of the tangent bundle of ${\cal T}$ is simply
$\pi^*{\cal O}(D)$.
This has a $deg(D)$-dimensional family of sections lifted from $\Sigma$;
it also has a tautological
section $\lambda$, which vanishes along the zero section in ${\cal O}(D)$.
  Each of these sections defines a Poisson structure
on ${\cal T}$; the one we will use is $\lambda$.

 In turn, each Poisson structure on the surface ${\cal T}$ induces a
Poisson
structure on moduli spaces of sheaves on ${\cal T}$  [Bo2,T]. 
The moduli space we consider is that of the sheaves $L$
defined above, which are
supported along the spectral curves. The tangent space to the moduli space
at $L$ is $Ext^1(L,L)$;
dually, the cotangent space is $Ext^1(L,L\otimes K_{\cal T})$. The Poisson 
structure can be thought 
of as a skew
map from the cotangent space to the tangent space; it is given here by 
multiplication by $\lambda$.
$$\hat\lambda: Ext^1(L,L\otimes K_{\cal T})\rightarrow Ext^1(L,L)\eqno
(3.8)$$
To compute the $Ext$-groups, one can first take a 
locally free  resolution $R$ of $L$, take 
the induced complex $Hom(R,L\otimes K_{\cal T})$, and then compute the
first hypercohomology group of this complex. 
We have already found a resolution; 
it is given by the sequence (3.4). Applying
$Hom$, and recalling that $K_{\cal T}= \pi^*{\cal O}(-D) $,
the  cotangent space will be the first hypercohomology of the
 complex 
$$ 
(\pi^* E)^*\otimes L\otimes \pi^*{\cal O}(-D)
{\buildrel{(g-\lambda\bbi)^*}
\over
{\longrightarrow}} (\pi^*E)^*\otimes L,\eqno (3.9)$$
and the tangent space will be   the first hypercohomology of
$$ (\pi^* E)^*\otimes L {\buildrel{(g-\lambda\bbi)^*}\over
{\longrightarrow}} (\pi^*E)^*\otimes L \otimes \pi^*{\cal O}(D).\eqno
(3.10)$$
The map between the two complexes 
is multiplication by the tautological section $\lambda$. 
Pushing this down to $\Sigma$ we obtain for the cotangent and  tangent
spaces the first hypercomology groups of 
$$ End(E)(-D) {\buildrel{-ad_g}\over
{\longrightarrow}} End(E),\quad  End(E)  {\buildrel{ad_g}\over
{\longrightarrow}} End(E)(D),\eqno (3.11)$$
respectively.
The map between  the complexes in (3.11)
is left multiplication by $g$.
The tangent and
cotangent spaces
 then fit into exact sequences
$$\matrix {0&\rightarrow& H^0(\Sigma, End(E))&\rightarrow & 
T^*_{(E,g)}{\cal M}
&\rightarrow &H^1(\Sigma, End(E)(-D))&\rightarrow&0\cr
 & & \downarrow&&\downarrow&&\downarrow&&\cr
0&\rightarrow& H^0(\Sigma, End(E)(D))&\rightarrow & T_{(E,g)}{\cal M} 
&\rightarrow &H^1(\Sigma, End(E))&\rightarrow&0.\cr}\eqno (3.12)$$
The vertical arrows are left multiplication by $g$.
Explicitly, in \v{C}ech terms with respect to a covering $U_\alpha$,
 the cocycles for the first 
hypercohomology group for a complex 
$A{\buildrel{\sigma}\over {\longrightarrow}} B$ 
are given by pairs $(a_{\alpha,\beta}, b_\alpha)$, where 
$a_{\alpha,\beta}$ is a 1-cocycle for $A$, and 
$b_\alpha$ a 0-cochain for $B$ satisfying 
$$\sigma (a_{\alpha,\beta}) - b_\alpha+b_\beta = 0\eqno (3.13)$$
on overlaps. The coboundaries in turn, are given by taking a cochain
$a_\alpha$ for $A$
and mapping it to $(a_\alpha-a_\beta, \sigma(a_\alpha))$.

Reduction by the $\bbc^*$-action (3.2) corresponds to fixing the 
top exterior power of $E$ and  taking the quotient by the 
$\bbc^*$-action.
A Zariski open subset of the moduli ${\cal M}$ consists of Higgs pairs
with
a stable vector bundle of degree $1$. 
Since such a stable vector bundle $E$ is unique up to tensoring by a line 
bundle, a component of 
the reduced moduli space is the projective space 
$\bbp H^0(\Sigma,End(E)(D))$, endowed with
a Poisson structure. We have the Casimir determinant 
morphism, from the generically invertible locus in 
$\bbp H^0(\Sigma,End(E)(D))$, to the linear system
$\bbp H^0(\Sigma,{\cal O}(N\cdot D))$. The generic fiber contains a 
maximal dimensional symplectic leaf. Fix a non-zero section $\delta$ of
$H^0(\Sigma,{\cal O}(N\cdot D))$. 
The locus $H^0(\Sigma,End(E)(D))_\delta$, of sections with determinant 
$\delta$, is a cyclic $N$-sheeted \'{e}tale covering of the symplectic 
leaf 
in $\bbp H^0(\Sigma,End(E)(D))$ determined by the 
zero divisor of $\delta$. 
We abuse notation and denote this symplectic cyclic cover by 
${\cal M}_{\rm red}(E,\delta)$, or ${\cal M}_{\rm red}$ for short.

Next we identify the tangent and cotangent spaces of ${\cal M}_{\rm red}$.
Denote by $End^{~0}(E)$ the subbundle of traceless endomorphisms. 
Let $End^g(E)$ be the subbundle of $End(E)$, which, 
away from the singularities of $g$, is the image of $End^{~0}(E)$ under
right multiplication by $g$ (left multiplication results with the same 
subbundle). $End^g(E)$ is the subsheaf of $End(E)$ of sections 
satisfying
$$
\{f \in End(E) \ : \ {\rm tr}(g^{-1}f)=0\}. \eqno (3.14) 
$$
If $g^{-1}$ is a nowhere vanishing holomorphic section of $End(E)(D')$, 
then it defines a line subbundle $L$ of $End(E)$ isomorphic to 
${\cal O}_\Sigma(-D')$. 
$End^g(E)$ is the subbundle $L^\perp$ orthogonal to $L$ with respect to
the 
trace pairing. It is isomorphic to the dual of the quotient
$End(E)/L$. Thus, ${\rm deg}(End^g(E))={\rm deg}(L)=-{\rm deg}(D')$. 
The tangent space of ${\cal M}_{\rm red}$ at $g$ is given by the first 
hypercohomology of the complex (in degrees $0$ and $1$)
$$
End^{~0}(E)  {\buildrel{ad_g}\over
{\longrightarrow}} End^g(E)(D).
$$
The cotangent space is given by the first 
hypercohomology of the dual complex (in degrees $0$ and $1$)
$$
End^g(E)^*(-D) {\buildrel{-ad_g^*}\over
{\longrightarrow}} End^{~0}(E)^*. 
$$
The Poisson structure is induced by a 
homomorphism $\Lambda$ from the cotangent complex to the tangent complex. 
In degree $1$, $\Lambda_1$ is the composition of the isomorphism 
$End^{~0}(E)^*\cong End^{~0}(E)$
with left multiplication by $g$ (which takes $End^{~0}(E)$ into 
$End^{~g}(E)(D)$).
It is simpler to describe the dual of the 
homomorphism $\Lambda_0$ in degree $0$. $\Lambda_0^*$ is the composition
of the isomorphism $End^{~0}(E)^*\cong End^{~0}(E)$ with {\it right} 
multiplication by $g$. The commutativity 
$ad_g\circ \Lambda_0 = - \Lambda_1\circ ad_g^*$
follows from that in the $Gl(N)$ case. (Note that the transpose of
right multiplication by $g$ is given by left multiplication by $g$ and
the transpose of $ad_g$ is $-ad_g$). 

Taking the first hypercohomologies and recalling that both
$H^0(\Sigma, End^{~0}(E))$ and $H^1(\Sigma, End^{~0}(E))$ vanish,
we find: 
$$\matrix {& & T^*_{(E,g)}{\cal M}_{\rm red}
&\simeq &H^1(\Sigma, End^{~g}(E)^*(-D))\cr
   &&\downarrow&& \cr
 H^0(\Sigma, End^{~g}(E)(D))&\simeq & T_{(E,g)}{\cal M}_{\rm red}  
&&\cr}\eqno
(3.15)$$
This procedure endows the Zariski open subset of $H^0(\Sigma,End(E)(D))$, 
of generically invertible sections, with a Poisson structure. 
As a homomorphism 
from $T^*_{(E,g)}H^0(\Sigma,End(E)(D))$ $ = H^1(End(E)(-D))$ to
$H^0(\Sigma,End(E)(D)) = T_{(E,g)}H^0(\Sigma,End(E)(D)),$ 
it factors through the homomorphisms (3.15).

{\sc Lemma} (3.16) \tensl 
If $D>0$ and $N>1$, the Poisson structure extends to 
the whole of $H^0(\Sigma,End(E)(D))$.
\tenrm

{\sc Proof:} One shows that the locus of non-invertible sections of 
$H^0(\Sigma,End(E)(D))$ has codimension $\geq 2$. 
Since $N>1$, it suffices to estimate the codimension in the subspace 
of traceless sections. 
Let $End^{~0}(E_{D})$ be its restriction to $D$. The evaluation
homomorphism
$End^{~0}(E)(D)\rightarrow End^{~0}(E_{D})$ is an isomorphism because
$H^1(End^{~0}(E))=0$. The determinant divisor in $\gsl_N$ is irreducible. 
If $D>0$ and $x$ is a point in $D$, we get an {\it irreducible} divisor in
$H^0(End^{~0}(E_{D}))$ of sections which are not invertible at $x$. 
It suffices  that one of those sections $\varphi$ is generically 
invertible. 
Indeed, a line bundle $L$ on a reduced and irreducible sectral 
curve passing through the zero point in the fiber over $x$ will give rise
to such a section $\varphi$. 
\hfill\za

{\sc Remarks.} 1) The left multiplication, appearing in the 
construction of the Poisson structure, corresponds to an embedding of the
Lie group $GL(N)$ in its Lie algebra. This embedding has been implicitly 
used  
when we described meromorphic elements of the loop group as Higgs fields, 
i.e., as meromorphic sections of a Lie algebra bundle. 

2) We could use, instead, right multiplication.
The resulting Poisson structures will 
be equal to the one coming from left multiplication.
Indeed, their difference $ad_g$ is a  homomorphism 
between  the complexes in (3.11), which is {\it homotopic to zero}. The
homotopy $h$, as a homomorphism of degree $-1$
between the complexes, is given by the identity from
$End(E)$ to $End(E)$. 

{\it 3.c. Comparing the Sklyanin and the Mukai brackets.}
We start with an element $c\in H^1(\Sigma, End^{~g}(E)(-D)) =
H^0(\Sigma, End^{~g}(E)^*(D))^*$.
 We choose an open cover $U_+, U_-$ compatible with $D$ a divisor 
 disjoint from the open disk $U_-$.
We can represent $c$ as a cocycle $c_{\pm}$ with respect to our cover.
Lifting to 
$Ext^1$, we have a class represented by $(c_{\pm}, \rho_+,\rho_-)$, with 
$gc_{\pm}-c_{\pm}g - \rho_+ + \rho_- = 0$ on $U_+\cap U_-$. Now note that 
$gc, cg\in  H^1(\Sigma, End(E))$ can be split as 
$$gc_{\pm}= \mu_+ -\mu_- + {1\over N}tr(gc_{\pm})\bbi,\quad cg_{\pm}=
\nu_+ -\nu_-
+ {1\over N}tr(cg_{\pm})\bbi,\eqno (3.17) $$
since $H^1(\Sigma, End^{~0}(E)) = 0$. 
The hypercohomology cocycle condition implies 
that one can choose $\mu_+,\mu_-,\nu_+,\nu_-$ to satisfy 
$\rho_+ =\mu_+ -\nu_+$ and $\rho_-= \mu_- -\nu_-$. 
With this, we can compute the explicit form of the
Poisson structure $\Lambda$
$$\eqalign {\Lambda: T^*{\cal M}_{\rm red} &\rightarrow T{\cal M}_{\rm
red}\cr
((c_{\pm}, \mu_+ -\nu_+,\mu_- -\nu_-)&\mapsto (gc_{\pm}, g\mu_+
-g\nu_+,g\mu_- -g\nu_-)}\eqno (3.18)$$ 
Now we modify the class on the right by the coboundary $-(\mu_+-\mu_-,
g\mu_+ - \mu_+g, g\mu_- - \mu_-g)$,
which rewrites the map (3.18) as:
$$
((c_{\pm}, \mu_+ -\nu_+,\mu_- -\nu_-) \mapsto (0, \mu_+g -g\nu_+,\mu_-g
-g\nu_-). \eqno (3.19)$$
 The cocycle condition on the right hand side of (3.19) tells us 
that $\mu_+g-g\nu_+ = \mu_-g-g\nu_-$, and so defines  a global section of 
$H^0(\Sigma,End^g(E)(D))$.

 Given two classes $c$ and $d$ in the cotangent space
$H^1(\Sigma,End^g(E)^*(-D))$ at $(E,g)$, the Poisson structure is 
given by 
$$<\Lambda(c),d> \ = \ <\mu_+g-g\nu_+,d_\pm> \ = \ 
<\mu_+ - Ad_g(\nu_+),gd_\pm>.
 \eqno (3.20)$$

Let us compute the Poisson bracket corresponding to this, on a pair of
functions
$f, h$ on $ H^0(\Sigma, End^{~g}(E) (D))$. The differentials $df, dh$ of
these functions at $g$ are naturally identified
with classes in $H^1(\Sigma, End^{~g}(E)^*(-D))$ via Serre's Duality
and the trace pairing.  We will need the following elementary Lemma.

{\sc Lemma (3.21) } 
\tensl 
Trivialize the tangent bundle of $Gl(N)$ via the inclusion 
$Gl(N)\subset gl(N)$. 
Let $\rho_g : GL(N) \rightarrow GL(N)$ denote the right multiplication by
$g$. 
Identify a one form $df$ on 
$GL(N)$ with a vector field $\phi$ via the above trivialization and the
trace multiplication pairing: 
$$<\xi,df> \ = \ tr(\xi\cdot \phi), \ \ \forall \ \  \xi\in gl(N).$$
Then the pull back of a $1$-forms $df$  by $\rho_g$ corresponds to left 
multiplication of $\phi$ by $g$.
\tenrm

{\sc Proof: }\ 
$<\xi,d\rho_g^*(df)> \ = \ <d\rho_{g}(\xi),df> \ = \ 
tr(\xi\cdot g \cdot \phi).
$

 For an infinitesimal variation $\dot g$ through $g$, 
$$< Df, \dot g g^{-1}> = <df,\dot g> = <D'f, g^{-1}\dot g>.$$
Thus, $(Df)(g)$ is identified with $d\rho_g^*(df)$. Using the
above Lemma, we can identify $Df$ with $g\cdot df$, 
and $D'f$ with $df\cdot g$, and similarly for $dh$. In particular,
 if we represent $df$ by a $1$-cocycle $c_\pm$ as above, 
then $P_+(Df)=P_+(gc_\pm)=\mu_+$. Similarly, we have:
$$P_\pm(Df) = \mu_\pm, \quad P_\pm(D'f) = \nu_\pm.\eqno(3.22)$$
Substituting this into the expression  (3.20) for the Poisson bracket 
gives 
$$\{f, h\} = <\Lambda(df) , dh> = <P_+(Df)-Ad_g(P_+(D'f)), Dh>.
 \ \ \   \eqno(3.23)$$

As this is the expression given above in (2.13) for the Sklyanin bracket,
we have:

{\sc Theorem (3.24) } \tensl The Mukai bracket and the Sklyanin bracket
coincide on the reduced symplectic leaf of $H^0(\Sigma, End(E) (D))$ 
consisting of endomorphisms with a fixed determinant $\delta$.
\tenrm

{\sc Remark:} 
1) There is a natural $\bbc^\times$-action on the surface ${\cal T}$,
and consequently on the moduli spaces ${\cal M}$ and ${\cal M}_{\rm red}$ 
of sheaves on ${\cal T}$. 
The Poisson structure we constructed is $\bbc^\times$-invariant
with respect to the natural $\bbc^\times$-action on 
$H^0({\cal M},{\buildrel{2}\over{\wedge}}T{\cal M})$. 
So is the Poisson structure we started with on the surface ${\cal T}$. The 
$\bbc^\times$-invariance is related to the quadratic nature 
of the Poisson structure. Indeed, Lemma (3.16) produced a Poisson 
structure 
on the vector space $V=H^0(\Sigma, End(E) (D))$, which must come from an 
element of $Sym^2(V^*)\otimes {\buildrel{2}\over{\wedge}}V$.

2) Polishchuk constructed a related quadratic Poisson structure
on the moduli space ${\cal N}$ of stable triples 
$(E_1,E_2,\phi:E_2\rightarrow E_1)$
(see [Po]). There is a natural morphism from our moduli space ${\cal M}$
of stable Higgs pairs to ${\cal N}$ 
(it involves taking the quotient my the $\bbc^\times$-action). 
The morphism is Poisson. 

\noindent
{\it 3.d. Birational symplectic isomorphisms with Hilbert schemes.} 
We can compute simple Darboux coordinates for the  Mukai symplectic form  
in the $Gl(N,\bbc)$-case. This will, incidentally, also show explicitly 
that we do have an integrable system, as well as characterise the 
symplectic leaves.  To do this, we construct different 
resolutions for $L$, to compute the 
Ext-groups of (3.8). Let us extend $L$ to a sheaf $L_U$ defined on 
an analytic neighbourhood $U$ of 
a smooth spectral curve. We then have the resolution, on $U$:
$$0\rightarrow L_U(-nD){\buildrel{det(g-\lambda\bbi)}\over
{\longrightarrow}}
L_U\rightarrow L\rightarrow 0\eqno (3.25)$$
taking duals, and tensoring with $L(-D)$, the cotangent space of our
moduli
will
be the hypercohomology, over $S$, of the sequence
$${\cal O}(-D)\rightarrow {\cal O}((n-1)D)\eqno (3.26)$$
and the tangent space, that of the sequence
$${\cal O} \rightarrow {\cal O}(nD)\eqno (3.27)$$
The maps in (3.26), (3.27)  induced by (3.25) are simply the zero map, and
so,
since $K_S = {\cal O}((n-1)D) $, the cotangent 
space splits as  
$$T^* {\cal M}\simeq   H^1(S,{\cal O}(-D))\oplus H^0 (S, K_S 
)\eqno(3.28)$$
and the tangent space as
$$T {\cal M}\simeq   H^1(S,{\cal O})\oplus H^0 (S, K_S(D))\eqno(3.29)$$
The Poisson structure $\Lambda: T^*{\cal M}\rightarrow T {\cal M}$ is
given, as above, by  multiplication by the tautological 
section $\lambda$ of ${\cal O}(D)$ on both summands. 

Let us define a subspace ${\cal M}^Z$ of ${\cal M}$ of pairs $(E,g)$ 
whose spectral curve intersects the zero section in ${\cal T}$ in a fixed 
divisor $Z$. ${\cal M}^Z$, by (2.14), is
a union of   symplectic leaves 
of the Poisson structure. Now identify $H^0 (S, K_S)$ as the subspace
 of $H^0 (S, K_S(D))$
of sections vanishing along the zero-section.  We have
$$T{\cal M}^Z\simeq   H^1(S,{\cal O})\oplus H^0 (S,
K_S)\eqno(3.30)$$
and dually
$$T^* {\cal M}^Z\simeq   H^1(S,{\cal O})\oplus H^0 (S,
K_S)\eqno(3.31)$$

Under our identifications, the Poisson structure $T^* {\cal
M}^Z\rightarrow
T {\cal M}^Z$ is simply the identity map. In other words, the
Poisson tensor
is the canonical one on the sum (3.30),  under the identifications we have
made. As the Serre pairing is non-degenerate, this shows, fairly
immediately, several important things:

{\sc Theorem (3.32)}: \tensl The space ${\cal M}^Z$ has an open set which 
is a symplectic leaf. The foliation on this leaf obtained by fixing the
spectral curve 
 is Lagrangian. The dimension of ${\cal M}^Z$ is twice the genus $g$ of
the 
 spectral curve, which is given by \tenrm
 $$g={N(N-1)d\over 2} +1$$
 The constant $d$ is the degree of the divisor $D$. The computation of 
 the genus is a simple application of the adjunction formula.

The canonical forms have  simple   Darboux coordinates. Indeed, the
bundle $E$
has a one dimensional space of  sections [At], and so, by (3.6), does $L$. 
In other words, the bundle $L$ is represented in a unique way as a divisor
$\sum_\mu p_\mu,\ p_\mu\in S$. Now we note that the space ${\cal T}$
admits
a symplectic form, unique up to scale, with a pole along the zero-divisor.
If $\lambda$ is a linear fibre coordinate on ${\cal T}$, and $\omega= dz$
is
 a one-form on the base elliptic curve $\Sigma$ (where $z$ is a standard
linear coordinate 
on the curve),  the symplectic form on $\Sigma$
is given by $\Omega_{\cal T}= {d\lambda\over\lambda}\wedge \pi^*\omega$.
The points $p_\mu$, which we will suppose distinct (as they are, 
generically), are given by a pair
of coordinates $( z_\mu,\lambda_\mu)$. These pairs not only determine the
line bundle but, generically, 
also the curve $S$, as it must pass through the points $p_\mu$.

{\sc Proposition } (3.33): \tensl The Mukai form on the symplectic leaves
can be written
as $\Omega = \sum_\mu  {d\lambda_\mu\over\lambda_\mu}\wedge dz_\mu$.\tenrm

{\sc Proof: } The proof follows verbatim that given, e.g., in [HK]. It is
mostly
a question of writing down the explicit form of the duality pairings.
\hfill\za

More invariantly, the proposition establishes a birational
symplectomorphism between the open
 symplectic leaves of $ {\cal M}^Z$  and the $g$-th
symmetric product
of a blow up $\widetilde {\cal T}$ of ${\cal T}$.
Indeed, along $ {\cal M}^Z$, the intersection $S\cap Z$ of the
spectral curves
with the zero-section is fixed. Let us blow up the
points
of $S\cap Z$, and call the resulting surface $\widetilde {\cal T}$. The
form 
$\Omega_{\cal T}$ lifts to a form $\widetilde {\Omega_{\cal T}}$ on
$\widetilde {\cal T}$, which is holomorphic away from the proper transform
of the zero section. 
Let $Hilb_g(\widetilde {\cal T})$ denote the Hilbert scheme 
of 0-dimensional length $g$ subschemes  
in $\widetilde {\cal T}$; this is a 
desingularisation of the symmetric product, and it is symplectic. 
Proposition (3.33) then becomes:

{\sc Proposition (3.34)}: \tensl On the  generic symplectic leaves 
 of the Mukai
bracket, the map which associates to a pair $(S,L)$ its divisor
$\sum_\mu p_\mu$ is a birational  symplectic map  between ${\cal L}$ and
$Hilb_g(\widetilde {\cal T})$. \tenrm

To deal with the Sklyanin bracket, we must reduce, both on the space of
sections of 
$End(E)(D)$ and on the Hilbert scheme. For the first, as we indicated, 
the reduction amounts to fixing the top exterior power of
$E$, 
and then quotienting by the action of $\bbc^*$ on the section $g$; 
equivalently, up to a finite cover, we fix the scale of 
the determinant; its zeroes are fixed on the symplectic leaf.
For the Hilbert scheme, the surface 
$\widetilde {\cal T}$ admits a $\bbc^*$-action along the fibers
of the projection $\widetilde {\cal T}\rightarrow \Sigma$.
This action is symplectic, and its moment map (with values in $\Sigma$)
is given by projection. The action extends to $Hilb_g(\widetilde {\cal
T})$;
 the moment map is then the sum in
$\Sigma$ of the  points $\pi(p_\mu)$. To reduce under this action, we must
fix the sum of the points, and then quotient by the $\bbc^*$ action. 
Note that, as in [HK], the sum of the points in $\Sigma$ is essentially
the divisor corresponding to the top exterior power of the push-down 
$E$ of the line bundle $L$.   Fixing the determinant of $g$, 
once one has its zeroes, results in a cyclic cover of 
the quotient by the $\bbc^*$ action. 
In short, the reductions by the $\bbc^*$ actions are compatible. We have:

{\sc Proposition (3.35)}: \tensl On the   symplectic leaves ${\cal L}$
 of the Sklyanin
bracket, the map which associates to a pair $(S,L)$ its divisor
$\sum_\mu p_\mu$ is a symplectic map  between ${\cal L}$ and
$Hilb_g(\widetilde {\cal T})//\bbc^*$.\tenrm

It is perhaps worth emphasizing that the above description is
 quite amenable to explicit calculation. Indeed, as we saw, 
using the projection $\pi:\bbc\rightarrow \Sigma$, elements of the 
symplectic leaf ${\cal L}$ can be described as matrix valued functions $M$
on 
$\bbc$ with poles at $\pi^{-1}(D)$ satisfying $M_i(z+\omega_i) =
I_iM(z)I_i^{-1}$;
these can be represented using theta-functions.
The points $(z_\mu,\lambda_\mu)$ can be computed as zeroes of the equation
$$(M(z)-\lambda\bbi)^{\rm adj}S=0,\eqno (3.36)$$
where adj denotes the matrix of cofactors, and  $S$
 is a column vector of functions representing the section of $E$. It can
be computed
explicitly using theta-functions, and the explicit formula is given in
[HK], 
section 4.

The coordinates $(z_\mu, \lambda_\mu)$ allow a simple linearisation of the
flows. 
Indeed, we note that fixing the Hamiltonians $H_1,...H_k$ fixes the
spectral curve, and
so determines $\lambda$ as a function of $z$: $\lambda = \lambda
(z,H_1,...,H_k)$.
Choosing a base point $z_0$ on $\Sigma$, we set
$$F(z_\mu, H_i) = \sum _\mu \int_{z_0}^{z_\mu} ln(\lambda(z, H_i))dz.
\eqno (3.37)$$
Since $\partial F/\partial z_\mu = ln(\lambda_\mu)$, the linearising
coordinates 
of the flows are 
given by 
$$Q_i = {\partial F\over \partial H_i} = \sum _\mu \int_{z_0}^{z_\mu}
\lambda^{-1}
{\partial \lambda\over \partial H_i}. \eqno (3.38)$$
One can show that these are sums of Abelian integrals.
\bigskip

{\bf 4. Rational nodal, or trigonometric case.} 

One can allow the elliptic
curve 
$\Sigma $ to degenerate, and obtain a rational nodal curve $\Sigma_0$
which is equivalent to 
the Riemann sphere $\bbp^1$ with two points $z=0, \infty$ identified. We
take
the bundle $ {\cal O}\oplus {\cal O}\oplus ... \oplus {\cal O}(1)$ 
of degree one on  $\bbp^1$, and identify the fibers over $0$, $\infty$ to
obtain a bundle $E$
on the rational nodal curve.
If one takes the transition matrix from $z\ne \infty$ to $z\ne 0$
$$T(z)= \pmatrix {0&1&0&\dots&0\cr 0&0&1&\dots&0\cr.&.&.& &.\cr.&.&.&
&.\cr
0&0&0&\dots&1\cr z^{-1}&0&0&\dots&0},\eqno(4.1)$$
the identification between the fibers can be taken to be the identity
matrix.
Alternately, we can pass to the universal cover $\bbc$ of $\bbc^*=
\bbp^1-\{0,\infty\}$
and use an automorphy factor representation, so that sections of $E$
are represented by vector functions satisfying $F(x+1) = I_1 F(x)$, and
suitable
boundary behaviour as $ix\rightarrow \pm \infty$. Endomorphisms again
become matrix valued functions
with $M(x+1) = I_1M(x) I_1^{-1}$. We refer to [HK], section 5.

Again these bundles are rigid, up to the top exterior power. There is
again a spectral curve
$S$, covering the 
curve $\Sigma_0$, and a line bundle $L$ on $S$, which can as above 
be represented by a divisor $\sum (z_\mu,\lambda_\mu)$. We can go through
the
proof of the identity of the reduced Mukai bracket with the Sklyanin
bracket,
 essentially verbatim. There is a  splitting of the loop group into the
sum of two
subgroups, one corresponding to sections on a neighbourhood
of $x=0$ (that is, $z=1$), and the other to sections of $End(E)$ on the
complement of $z=1$.
Again, the Mukai symplectic form on the leaves has the form $\sum_\mu 
{d\lambda_\mu\over \lambda_\mu}
\wedge dz_\mu$; the reduction to the Sklyanin form amounts to fixing the 
determinant of the curve, and fixing the product of the $z_\mu$. The
formula 
for computing the $(z_\mu,\lambda_\mu)$ are similar.

\bigskip

{\bf 5. Rational case.}

 While this case has already been computed
explicitly in [Sc],
the proof given above adapts in a straightforward way to cover this case,
too.
Our bundle $E$, now, is simply the trivial rank $N$ bundle over $\bbp_1$.
The bundle is, indeed,
rigid; however, $H^0(\bbp^1, End(E))\ne 0$, and so there is no unique
splitting of the
sections of $End(E)$ over the punctured disk. The groups $H^0(\bbp^1,
End(E)(-1))$ and 
$H^1(\bbp^1, End(E)(-1))$ are zero, however, and this gives a
decomposition
of 
sections of $End(E)$ over the punctured disk into a direct sum of
\item {-} the subalgebra of sections of $End(E)$  over the 
disk which vanish at the origin, and 
\item{-} the subalgebra of sections of $End(E)$ which are defined on the
complement of
the origin.

The spectral curves of elements $g$ of $H^0(\Sigma, End(E) (D))$
lie in the total space ${\cal T}$ of the line bundle 
${\cal O}(D)$ over $\bbp^1$. Let $z$ be the standard coordinate on
$\bbp^1$, and
let $\lambda$ be a linear coordinate along the fibers of ${\cal T}$. 
 The symplectic leaves lying in $H^0(\Sigma, End(E) (D))$ correspond to 
sections with a fixed determinant, as well as spectrum, 
which is fixed to order two over the point at infinity in $\bbp^1$. 
In other words, the spectral curves
have fixed intersection with  the zero-section $\lambda = 0$, 
as well as with $(z^{-2}) = 0 $. This foliation by symplectic leaves 
corresponds to the choice of a Poisson structure on ${\cal T}$, whose
divisor is precisely the zero-section $\lambda=0$ and twice the fiber
over $z=\infty$.

The space $H^0(\Sigma, End(E) (D))$ is acted on by  
$PGL(N,\bbc)$ via the adjoint action of the group of automorphisms of the 
trivial bundle. We can take the Poisson quotient, to obtain a reduced 
space $H^0(\Sigma, End(E) (D))/PGL(N,\bbc)$. 
We have, in a fashion analoguous to what is given above:

{\sc Proposition (5.1) } \tensl The Mukai Poisson structure and the
reduced
Sklyanin
structure coincide. If elements $g$ correspond to a line bundle $L$ over
the spectral
curve, represented by a divisor $\sum_\mu (z_\mu,\lambda_\mu)$, the
symplectic form on the leaves is 
$\sum_\mu {d\lambda_\mu\over \lambda_\mu}\wedge dz_\mu$.\tenrm
 
\bigskip
{\bf 6. Higher genus}. 

One can ask how the above extends to higher genus base curves.
One still, of course, has a space of pairs $(E,g)$, consisting of rank $N$
bundles $E$ and sections $g$ of $H^0(\Sigma, End(E)(D))$. If $D$ is the 
sum of a canonical divisor and an effective divisor, then 
there is a Poisson structure on this
space. The Poisson structure corresponds to the generalised Hitchin 
systems. 
Following the procedure of Mukai, it corresponds to a Poisson structure on 
${\cal O}(D)$ which is constant
 along the fibers of the projection $\pi$ to the base curve $\Sigma$. The
Sklyanin systems, on the other
hand, correspond to Poisson structures which are linear along the fibers.
These
only exist if the genus is at most one; if the genus is greater, one only
has meromorphic
Poisson structures, of the form $\lambda{\partial \over \partial
\lambda}\wedge
\pi^*\omega^{-1}$, where $\omega$ is a holomorphic form on $\Sigma$, and 
$\lambda $is a coordinate along the fiber. These forms
correspond to 
degenerate symplectic forms on the Jacobian fibration $(S, L)\rightarrow
S$.
The form is null on certain directions in the fibers of the Jacobian
fibration:
if $Z$ denotes the zero locus of $\omega$, the  null direction in the
Jacobian corresponds
to the coboundary $\delta(H^0(S\cap\pi^{-1}(Z), {\cal O}(\pi^{-1}(Z))))$
in the exact sequence
$$...\rightarrow H^0(S\cap\pi^{-1}(Z), {\cal O}(\pi^{-1}(Z)))\rightarrow
H^1(S,{\cal O})
\rightarrow  H^1(S,{\cal O}(\pi^{-1}(Z))).$$
In any case, we can see that the Poisson geometry of rational surfaces
suggests quite strongly that
there is no nice Poisson extension of the Sklyanin bracket to arbitrary
base curves.

\bigskip

{\bf Bibliography}
\medskip

\parskip = 0 pt
\item{[At]}
M. Atiyah,
{\it Vector bundles over an elliptic curve},
Proc. Lond. Math. Soc {\bf 7} (1957), 414--452.

 \item{[AHH]}
M.R. Adams,J. Harnad and J. Hurtubise, 
{\it Darboux coordinates and Liouville-Arnold integration in loop
algebras}, 
Comm. Math. Phys. {\bf 155} (1993), no.~2, 385--413.

\item{[AHP]} 
M.R. Adams, J. Harnad and E. Previato, 
{\it Isospectral Hamiltonian flows in finite and infinite dimensions {\rm
I}.
Generalised Moser systems and moment maps into loop algebras}, 
Comm. Math. Phys. {\bf 117} (1988), no.~3, 451--500.

\item{[AvM]}
M. Adler and P. van Moerbeke, 
{\it Completely integrable systems, Euclidean Lie algebras, and curves}, 
Adv. in Math. {\bf 38} (1980), no.~3, 267-317; 
{\it Linearization of Hamiltonian systems, Jacobi varieties and
representation theory}, ibid. {\bf 38} (1980), no.~3, 318--379.

\item{[BD]} A. A. Belavin and  V. G. Drinfeld:
{\it Solutions of the classical Yang-Baxter equations for simple Lie 
algebras.\/}
Funct. Anal. and its appl., {\bf 16 } (1982), 159-180

\item{[Bo1]}
F. Bottacin, 
{\it Symplectic geometry on moduli spaces of stable pairs,}
Ann. Sci. Ecole Norm. Sup. (4) {\bf 28} (1995), no. 4, 391-433.

\item{[Bo2]}  F. Bottacin, 
{\it Poisson structures on moduli spaces of sheaves over Poisson
surfaces}, Invent. Math. {\bf 121} (1995), no. 2, 421-436. 


\item{[EV]} P. Etingof and A.Varchenko, {\it Geometry and classification 
of solutions to the classical dynamical Yang-Baxter equation}
Commun. Math. Phys. 192 (1998), no. 1, 77--120.

\item{[F]} G. Felder, {\it Conformal field theory and integrable systems 
associated to elliptic curves}, Proceedings of the ICM,
Vol. 1, 2 (Zürich, 1994), 1247--1255, Birkhäuser, Basel, 1995. 
 
\item{[FT]}
L.D. Faddeev and L.A. Takhtajan,
{\it Hamiltonian methods in the theory of solitons},
eds., Springer-Verlag, Berlin, 1987.

\item{[FM]}
R. Friedman and J. W. Morgan,
{\it Holomorphic principal bundles over elliptic curves},
math.AG/9811130

\item {[H]}
J. Hurtubise,
{\it Integrable systems and algebraic surfaces},  
Duke Math. J. {\bf 83} (1996), no.~1, 19--50.

\item{[HaH]} 
J. Harnad and J. Hurtubise, 
{\it Generalised tops and moment maps into loop algebras}, 
J. Math. Phys. {\bf 37} (1991), no.~7, 1780--1787.

\item{[Hi1]} 
N.J. Hitchin, 
{\it The self-duality equations on a Riemann surface}, 
Proc. London Math. Soc. (3) {\bf 55} (1987), no.~1, 59--126. 

\item{[Hi2]} 
N.J. Hitchin,  
{\it Stable bundles and integrable systems}, 
Duke Math. J. {\bf 54} (1987), no.~1, 91--114.


\item {[HM]}
J. Hurtubise and E. Markman, 
{\it The elliptic Sklyanin bracket for arbitrary reductive groups}, in
preparation.

\item{[M]}
E. Markman,   
{\it Spectral curves and integrable systems}, 
Compositio Math. {\bf 93} (1994), 255-290.

\item{[Mo]}J.  Moser,  ``Geometry of Quadrics and Spectral Theory", 
{\it The Chern Symposium, Berkeley, June 1979}, 147-188, Springer, New
York, (1980).

\item{[Mu]} S. Mukai,
{\it Symplectic structure of the moduli
space of sheaves on an abelian or K3 surface},
Invent. math. 77 (1984) 101-116

\item {[NV]} S.P.~Novikov and A.P. ~Veselov, {\it Poisson brackets and 
complex tori} Proc. of the Steklov Inst. of Math., {\bf 165} (1984),
53-65.

\item {[Po]} A. Polishchuk, 
{\it Poisson structures and birational morphisms associated with
bundles on elliptic curves}, Internat. Math. Res. Notices (1998), no. 
{\bf 13}, 683--703. 

\item{[RS1]} 
A.G. Reiman and M.A. Semenov-Tian-Shansky, 
{\it Reduction of Hamiltonian systems, affine Lie algebras and lax
equations {\rm I, II}}, 
Invent. Math. {\bf 54} (1979), no.~1,  81--100; ibid. {\bf 63} (1981),
no.~3, 423--432.

\item{[RS2]} 
A.G. Reiman and M.A. Semenov-Tian-Shansky, 
{\it Integrable Systems II}, chap.2, in ``Dynamical Systems VII'', 
Encyclopaedia of Mathematical Sciences, vol 16., 
V.I. Arnold and S.P.Novikov, eds., Springer-Verlag, Berlin, 1994.

\item {[S]} M.A. Semenov-Tian-Shansky {\it Dressing transformations and
Poisson group actions}
Publ. Res. Inst. Math Sci. {\bf 21} (1985), 1237-1260.

\item{[Sc]}
D.R.D. Scott, 
{\it Classical functional Bethe ansatz for $SL(N)$:separation of
variables for the magnetic chain}, 
J. Math. Phys. {\bf 35}, 5831-5843 (1994)

\item{[Sk1]}
E.K. Sklyanin,
{\it On the complete integrability  of the Landau-Lifschitz equation},
LOMI preprint E-3-79, (1979).

\item{[Sk2]}
E.K. Sklyanin,
{\it Poisson structure of a periodic classical $XYZ$-chain},
J. Sov. Math. ,{\bf  46} (1989) 1664-1683.

\item{[T]}  A. N. Tyurin,
{\it Symplectic structures on the varieties of moduli of
vector bundles on algebraic surfaces with $p_{g} > 0$},
Math. USSR Izvestiya Vol. {\bf 33}  No. 1 (1989) 139-177.

\ninerm 

\bigskip 
 
\line{   J. C. Hurtubise\hfil E. Markman} 
\line{Centre de Recherches Math\'ematiques\hfil Department of Mathematics}
\line{Universit\'e de Montr\'eal \hfil   University of Massachusetts}
\line {and Department of Mathematics\hfil Amherst} 
\line{McGill University\hfil  
email: markman@math.umass.edu} 
\line{email: hurtubis@crm.umontreal.ca \hfil }

\end